\documentclass[preprint,12pt]{elsarticle} 
\usepackage[english]{babel}
 
\usepackage{enumerate}           
\usepackage{amsmath}   
\usepackage{amssymb}   
\usepackage{graphicx}      
\usepackage{epsfig}            
     
\usepackage{color}       
\usepackage{caption}  
\usepackage{float}
\usepackage{tabulary}

\journal{...}  
\begin{document} 

\makeatletter
\newcommand{\norm}[1]{\left\lVert#1\right\rVert}
\makeatother
  

\begin{abstract}
In this brief note, we establish a novel criterion for robustness of global asymptotic stability of zero solution of LTV system $\dot x=A(t)x$ in the presence of possibly unbounded perturbations (external disturbances). To prove the result, logarithmic norm will be used under which the stability becomes a topological notion depending on the chosen vector norm in the state-space $\mathbb{R}^n.$
\end{abstract} 

\title{Criterion for robustness of global asymptotic stability to perturbations of linear time-varying systems}
\author{R.~Vrabel}
\ead{robert.vrabel@stuba.sk}
\address{Slovak University of Technology in Bratislava, Institute of Applied Informatics, Automation and Mechatronics,  Bottova 25,  917 01 Trnava,   Slovakia}

\newtheorem{thm}{Theorem}
\newtheorem{lem}[thm]{Lemma}
\newtheorem{defi}[thm]{Definition}
\newdefinition{rmk}{Remark}
\newdefinition{ex}{Example}
\newproof{pf}{Proof}

\pagestyle{headings}

\maketitle

\section[Introduction]{Introduction}

Let us consider $x=0$ being asymptotically stable equilibrium point of the linear time-varying (LTV) system $\dot x=A(t)x,$ $t\geq0.$ What can we say about the asymptotic stability of its perturbation $\dot x=A(t)x+w(x,t)?$ This question represents one of the fundamental problems in the area of  {\it robustness of stability to external perturbations (disturbances)} (cf. \cite{Goebel}) and robustness of the systems in general, and so the effect of (known or unknown) perturbations on the solutions of nominal system as a potential source of instability attracts the attention and interest of scientific community for a long time in the various contexts.  For example, robustness to the external perturbations is of the utmost importance for the robotic manipulators deployed in uncontrolled environments. Robustness \& resilience are often thought of in terms of a system's capacity to maintain its functionality (stability) in the face of external disturbances. A comprehensive overview of the most significant results on robust control theory and its history is presented in \cite{Petersen}.

In this paper, we specifically prove that asymptotically stable zero solution of unperturbed system remains "attractive" for a wide class of perturbations $w(x,t)$ in the sense of convergence of all solutions of perturbed system to the origin $x=0$ as $t\to\infty$ provided the perturbing term satisfies some growth constraints. Sufficient conditions ensuring this property of the systems, possibly even for (monotonically) {\bf unbounded} perturbations, are summarized in Theorem~\ref{thm:linear} and simulation experiment in Example~\ref{example2} confirms theoretical conclusions. 

This, maybe a little unexpected, behavior of the LTV systems follows from 
\begin{itemize}
\item[1)]  the complexity and richness phenomena of dynamics of unperturbed LTV systems, especially with unbounded $A(t),$ compared to the unperturbed LTI systems where the solutions are restricted to the subsets of linear subspace spanned by the narrow class of functions; 
\item[2)] the properties of convolution integral (representing system's response to external disturbances) in  Lagrange's variation of constants formula and its estimate in the terms of logarithmic norm, mixing perturbing term $w(x,t)$ with fundamental matrix solution $\Phi(t)$ (and its inverse) of $\dot x=A(t)x.$ 
\end{itemize}

As has been shown in \cite{Khalil}, if the origin $x=0$ is an exponentially stable equilibrium point of unperturbed system and the perturbation term $w$ satisfies
\begin{equation}\label{bound} 
\norm{w(x,t)}\leq\rho(t)\norm{x}+\xi(t),\ \forall\norm{x}<r,\ \forall t\geq0 
\end{equation}
where $\rho,\xi: [0,\infty)\rightarrow[0,\infty)$ are continuous, $\int_0^{\infty}\rho(\tau)d\tau<\infty$ and $\xi$ is bounded, then for $\xi\equiv0,$ the origin is an exponentially stable equilibrium point of perturbed system and the solutions of perturbed system are ultimately bounded in the opposite case, that is, if $\xi$ is not identically zero. The general framework for our considerations and analyses here is that we will not assume {\it a priori} that $w(x,t)$ satisfies the inequality constraint of the form (\ref{bound}) and therefore the now classic results of  Khalil \cite{Khalil} based on the Lyapunov's converse theorem, Coddington \& Levinson \cite{Coddington_Levinson}, Hartman \cite{Hartman} both based on the state-space model representation are not applicable here in general.

Since then it has not made any further substantial progress in the research of ability of the dynamical systems and among them LTV ones to absorb the effect of external disturbances and maintain the stability of ongoing processes.   
This is all the more surprising as the stability analysis for time-varying linear systems is of constant interest in the control community. One reason is the growing importance of adaptive controllers for which underlying closed-loop adaptive system is time-varying and linear \cite{Ioannou}, \cite{Shamma}. The second one is that LTV systems naturally arise when one linearizes nonlinear systems around a non-constant nominal trajectory. In contrast the linear time-invariant (LTI) cases which have been thoroughly understood in the analysis and synthesis, many properties of the LTV systems are still not clear and not resolved. Here, the system (robust) stability analysis can serve as an appropriate example. 

There are some papers providing some sufficient conditions for exponential, see e.~g., \cite{Ilchmann, Wu, Zhou}) or/and asymptotic stability \cite{Safavi, Wang} of LTV systems  but none of those deal with the robustness of systems' stability to external perturbations and, moreover,  the developed techniques are not directly applicable for perturbed systems.

\subsection{Notation and assumptions}

Throughout this paper we assume that $x$ and $w$ are $n-$dimensional column vectors and $A(t),$ $t\geq0$ is a square matrix of the same dimension.
We will always assume that $A(\cdot):\,[0,\infty)\to\mathbb{R}^{n\times n}$ is a continuous matrix function and $w$ is continuous in $(x,t)$ for $\norm{x}<\infty$ and $0\leq t<\infty.$  Notice that $x=0$ may not be solution of perturbed system.

We will derive the results for unspecified vector norm, $\norm{\cdot}.$ For the matrices, as an operator norm is always used the induced norm, $\norm{A}=\max\limits_{\norm{x}=1}\norm{Ax}$. For example, the Frobenius norm of the matrix is not induced norm because $\norm{I_n}=1$ for any induced norm, but $\norm{I_n}_F=\sqrt{n}.$  We use for both vector norm and matrix operator norm the same notation but it will always be clear from the context that norm is just being used. In particular cases we will consider the three most common vector norm: 
\begin{equation}\label{norms}
\norm{x}_1=\sum\limits_{i=1}^{n}|x_i|, \qquad \norm{x}_2=\sqrt{\sum\limits_{i=1}^{n}x_i^2}, \qquad \norm{x}_{\infty}=\max_{1 \leq i \leq n}{|x_i|}.
\end{equation}
We denote by $\mu[A(t)],$ $t\geq0,$  the logarithmic norm of matrix $A(t).$ The classical definition is
\[
\mu[A(t)]\triangleq\lim\limits_{h\to 0^+}\frac{\norm{I_n+hA(t)}-1}{h},
\]
where $I_n$ denotes the identity on $\mathbb{R}^n,$ \cite{Dekker_Verwer, Soderlind1, Soderlind2}.  To distinguish between the logarithmic norms for vector norms defined in (\ref{norms}) we use here the notation with the subscript $1,2$ and $\infty$.

Because it is assumed that $A(t)$ is continuous, the function $\mu[A(t)]$ is a continuous function of $t$ by virtue of the inequality $\vert\mu[A(t_1)]-\mu[A(t_2)]\vert\leq\norm{A(t_1)-A(t_2)},$ $t_1,t_2\geq0,$ see e.~g.~\cite{Coppel}.  While the matrix norm $\norm{A}$ is always positive if $A\neq0,$  the logarithmic norm $\mu (A)$ may also take negative values, e.~g. for the Euclidean vector norm $\norm{\cdot}_2$ and when $A$ is negative definite because $\frac12(A+A^T)$ is also negative definite, \cite[Corollary~14.2.7.]{Harville} and Table~\ref{table:norms}. Therefore, the logarithmic norm does not satisfy the axioms of a norm.

The superscript 'T' denotes transposition,  the number $\lambda_{\max}(M)$ in Table~\ref{table:norms} and elsewhere indicates the maximum eigenvalue of matrix $M.$

\centerline{}

\subsection{Preliminary results}

In this subsection we summarize important and interesting properties of the logarithmic norm useful for stability analysis of dynamical systems.

In Table~\ref{table:norms}, the values of matrix norms and corresponding logarithmic norms for the norms (\ref{norms}) are calculated, see e.~g. \cite[p.~54]{Afanasiev}, \cite[p.~33]{Desoer_Vidyasagar}.
\begin{table}[H]
\caption{Logarithmic norms for the vector norms $\norm{\cdot}_i,$ $i=1,2,\infty$}
\centering 
\scalebox{0.9}{
\begin{tabulary}{\linewidth}{LLL}
\hline\hline                        
Norm of vector ($\norm{x}_i$) & Induced norm of matrix $A$ ($\norm{A}_i$) & Logarithmic norm ($\mu_i[A]$)  \\  [0.5ex]
\hline               
$\norm{x}_1=\sum\limits_{i=1}^{n}|x_i|$ & $\norm{A}_1=\max\limits_{1 \leq j \leq n}\sum\limits_{i=1}^{n}|a_{ij}|$ \newline (column sum)  & $\mu_1[A]=\max\limits_{1 \leq j \leq n}\big(a_{jj}+\sum\limits_{i\neq j}|a_{ij}|\big)$   \\
  \\
$\norm{x}_2=\sqrt{\sum\limits_{i=1}^{n}x_i^2}$ & $\norm{A}_2=\sqrt{\lambda_{\max}(A^TA)}$ & $\mu_2[A]=\frac12\lambda_{\max}\left({A+A^T}\right)$   \\
  \\
$\norm{x}_{\infty}=\max\limits_{1 \leq i \leq n}{|x_i|}$ & $\norm{A}_{\infty}=\max\limits_{1 \leq i \leq n}\sum\limits_{j=1}^{n}|a_{ij}|$\newline  (row sum) & $\mu_{\infty}[A]=\max\limits_{1 \leq i \leq n}\big(a_{ii}+\sum\limits_{j\neq i}|a_{ij}|\big)$   \\  [1ex]   
\hline
\end{tabulary}
}
\label{table:norms}
\end{table}

\centerline{}

The fundamental advantage of approach based on the use of logarithmic norm is the fact that to estimate the norm of fundamental matrix solution $\Phi(t),$ $t\geq0$ for unperturbed system $\dot x=A(t)x$ and state-transition matrix $\Phi_{\mathrm{trans}}(t,\tau)\triangleq\Phi(t)\Phi^{-1}(\tau)$  we do not need to know its solutions explicitly. 

\centerline{}

\begin{lem}[\cite{Desoer_Vidyasagar, Soderlind1, Soderlind2}]\label{lognorm_properties}
\begin{itemize}
\item[]
\item[p1)] $\max\limits_{1\leq i\leq n}\Re(\lambda_i)\leq\mu[A],$ where $\lambda_1,\dots,\lambda_n$ are  the eigenvalues of $A,$ and $\Re(\lambda)$ is a real part of $\lambda;$
\item[p2)] $-\norm{A}\leq\mu[A]\leq\norm{A};$
\item[p3)] Let $\Phi(t),$ $t\geq0$ is a fundamental matrix solution for $\dot x=A(t)x.$  Then 
\[
\norm{\Phi(t)\Phi^{-1}(\tau)}\leq e^{\int\limits_{\tau}^t \mu[A(s)]ds}
\]
for all $t\geq\tau\geq0.$
\end{itemize}
\end{lem}
This lemma, Item~p3, allows to estimate the norm of state-transition matrix $\Phi_{\mathrm{trans}}(t,\tau)=\Phi(t)\Phi^{-1}(\tau)$ without knowing the fundamental matrix solution, purely on the basis of the matrix $A(t)$ entries, which can be especially important if $A$ is a non-constant matrix.

The following example demonstrates that the value $\mu[A]$ may depend on the used vector norm.
\begin{ex} \cite[p.~56]{Afanasiev} \label{example_different_norms}
\begin{itemize}
\item[a)]  
$
A_1=\begin{bmatrix}
-11 & 10 \\
 2 & -3 
\end{bmatrix}\Rightarrow 
$
$\mu_1[A_1]= 7,$ $\mu_2[A_1]=0.2111$ and $\mu_\infty[A_1]=-1;$
\item[b)]  
$
A_2=\begin{bmatrix}
-11 & 2 \\
 10 & -3 
\end{bmatrix}\Rightarrow
$
$\mu_1[A_2]= -1,$ $\mu_2[A_2]=0.2111$ and $\mu_\infty[A_2]=7;$
\item[c)]  
$
A_3=\begin{bmatrix}
-1 & 3 \\
 -3 & -2 
\end{bmatrix}\Rightarrow
$
$\mu_1[A_3]= 2,$ $\mu_2[A_3]=-1$ and $\mu_\infty[A_3]=2.$
\end{itemize}
Thus, we can verify whether the system $\dot x=A_ix, $ $i=1,2,3$ is stable or not by means of the vector norm with negative value of $\mu[A_i].$ In general, we obtain such logarithmic norm for any Hurwitz matrix $A$ \cite[p.~135]{Khalil}  using a vector norm $\norm{x}_H=\sqrt{x^THx},$ where the symmetric positive definite matrix $H$ satisfies the Lyapunov equation $A^TH + HA = -2I_n.$ The corresponding logarithmic norm $\mu_H[A]=-1/\lambda_{\max}(H),$ see Lemma~2.3 in \cite{Hu_Liu}. It turns out that the stability analysis based on the logarithmic norm becomes a topological notion unlike the spectrum of matrices which is topological invariant.
\end{ex} 
\section[robust global stability]{Robust global asymptotic stability of LTV systems}
Let us define the following classes of functions \cite{Strauss_Yorke}.
\begin{defi} 
Let $h: [0,\infty)\to\mathbb{R}^n$ be continuous. Define

\[
{\cal{V}}\triangleq\left\{ h:\ \norm{h(t)}\to0\quad \mathrm{as} \quad t\to\infty  \right\},
\]
\[
{\cal{AD}}\triangleq\left\{h: \int\limits_t^{t+1}\norm{h(s)}ds\to 0\quad \mathrm{as} \quad t\to\infty \right\},
\]
\[
{\cal{D}}\triangleq\left\{h: \sup\limits_{0\leq\eta \leq1}\norm{\int\limits_t^{t+\eta}h(s)ds}\to 0\quad \mathrm{as} \quad t\to\infty \right\}.
\]
\end{defi}
\begin{lem}\label{subset}
$
\cal{V}\subsetneq\cal{AD}\subsetneq\cal{D}. 
$
\end{lem} 
\begin{pf}
Assume that $\norm{h(t)}\to0$ as $t\to\infty.$ Using the monotonicity of the function $H(t)=\int\limits_0^t \norm{h(s)}ds$ and Lagrange's Mean Value Theorem we obtain the chain of inequalities from which immediately follows the claim of lemma,
\[
\sup\limits_{0\leq\eta \leq1}\norm{\int\limits_t^{t+\eta}h(s)ds}\leq \sup\limits_{0\leq\eta \leq1}\int\limits_t^{t+\eta}\norm{h(s)}ds = \int\limits_t^{t+1}\norm{h(s)}ds=\norm{h(\xi)},
\]
where $\xi\in(t,t+1),$ $\xi\to\infty$ as $t\to \infty.$ Now we need only to show that ${\cal{V}}\neq{\cal{AD}}\neq{\cal{D}}.$
\begin{itemize}
\item[i)] Let $h(t)=\left[\sin (e^t),\cos (e^t),0,\dots,0\right]^T.$ Then $h\in{\cal{D}};$ in fact, for any $\eta\geq0,$ 
\[
\left\vert\int\limits_{t}^{t+\eta}\sin (e^\tau) d\tau\right\vert\leq 2e^{-t}(1+e^{-\eta})\leq 4e^{-t}.
\]
The same inequality holds for the second component of $h,$ and thus
\[
\norm{\int_t^{t+\eta} h(s)ds}_2\leq\sqrt{32}e^{-t}\rightarrow 0\quad \mathrm{as}\quad t\to\infty,
\] 
but 
\[
\int_t^{t+1}\norm{h(s)}_2 ds=1 \not\to 0\quad \mathrm{as}\quad t\to\infty,
\]
that is, $h\in{\cal{D}},$  but $h\notin{\cal{AD}}.$ Recall that the use of the Euclidean norm in the last two steps does not impair the generality of analysis because all norms on finite-dimensional vector space are equivalent. Specifically, there exists a pair of real numbers $0<C_1\leq C_2$ such that, for all $x\in\mathbb{R}^n,$ the following inequality holds \cite{Johnson}:
\[
C_1\norm{x}_a\leq\norm{x}_b\leq C_2\norm{x}_a.
\]
In particular,
\[
\norm{x}_{\infty}\leq\norm{x}_{2}\leq\norm{x}_{1}\leq\sqrt{n}\norm{x}_{2}\leq n\norm{x}_{\infty}.
\]
\item[ii)] Let $h(t)=\left[h_1(t),0,\dots,0\right]^T,$ where $h_1(t)=\sum\limits_{n=1}^{\infty}h_{1,n}(t)$  is  "needle-like" function with $h_{1,n}$ defined as follows:
\[
h_{1,n}(t)=\left\{ 
\begin{array}{l}
    2n(t-n+1) \quad \mathrm{for}\ t\in\left[n-1,n-1+\frac{1}{2n}\right)\\
    2(-nt+n^2-n+1) \quad \mathrm{for}\ t\in\left[n-1+\frac{1}{2n},n-1+\frac{1}{n}\right)\\
    0\qquad\qquad\qquad\mathrm{elsewhere}
\end{array} \right. 
.\] 
Then 
\[
\int_t^{t+1}h_1(s) ds\leq\frac{1}{2n}\to 0\quad \mathrm{as}\quad t\to\infty,
\]
but $\norm{h(t)}\not\to0$ as $t\to\infty$ because for all $n\in\mathbb{N},$ the value of $h_1$ at the points $n-1+1/2n,$ $n\in\mathbb{N}$ is equal to $1;$ $h\in{\cal{AD}},$ but $h\notin{\cal{V}}.$
\end{itemize} 

Graphs of the first components of vector functions $h$ from the example above are shown in Figure~\ref{functions_h_1}.
\begin{figure}[!ht] 
\captionsetup{singlelinecheck=off} 
   \centerline{
    \hbox{
     \psfig{file=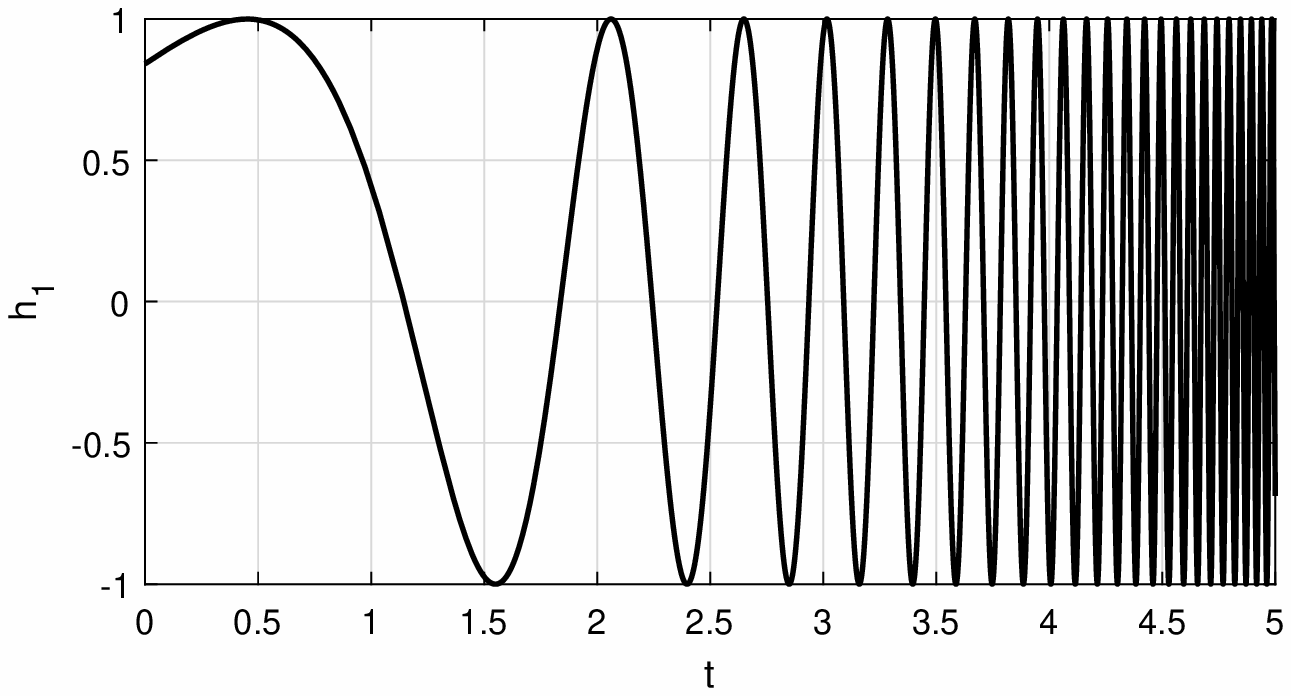,width=5.0cm, clip=}
     \hspace{1.cm}
     \psfig{file=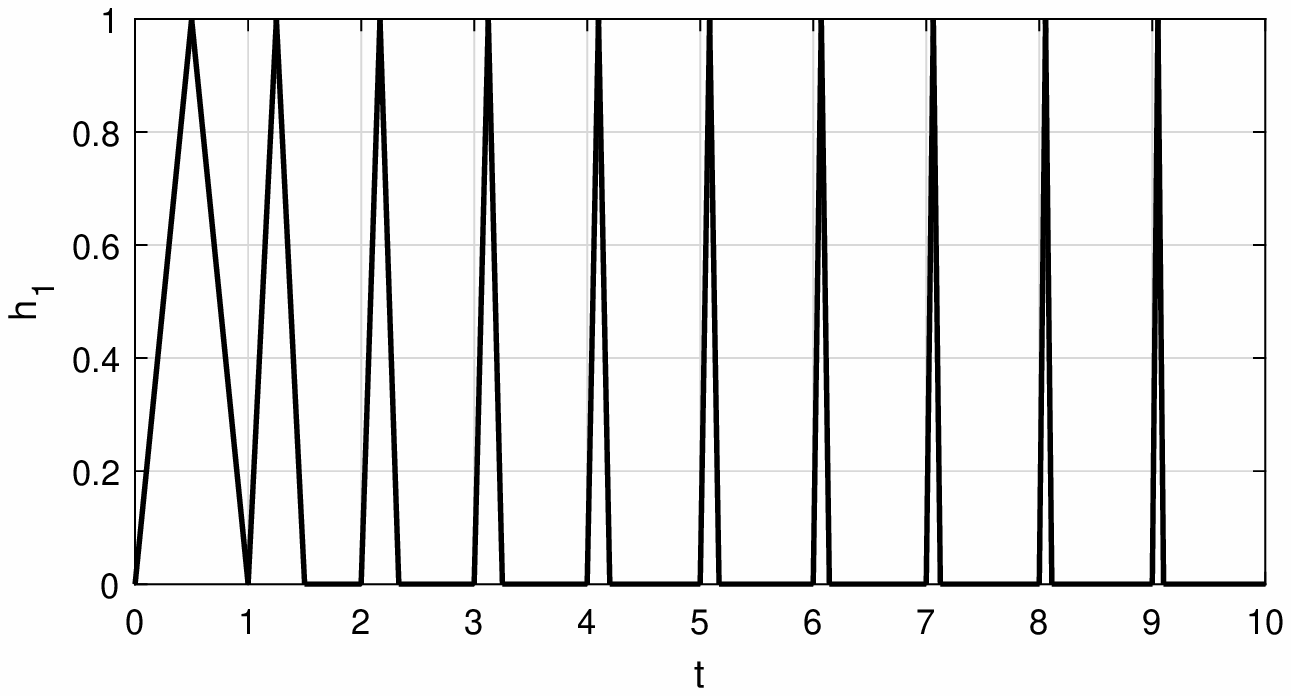,width=5.0cm,clip=}
    }
   }
   
\caption{The functions $h_1(t)$ from the proof of Lemma~\ref{subset}, i) on the left and ii) on the right.}
\label{functions_h_1}
\end{figure} 
This completes the proof of Lemma~\ref{subset}.
\end{pf}
We have the following result regarding asymptotic behavior of perturbed linear systems.
\begin{thm}\label{LTV} Let $A(\cdot):\,[0,\infty)\to\mathbb{R}^{n\times n}$ is a continuous matrix function. 
Let $x=0$ is globally uniformly asymptotically stable equilibrium point of the LTV system $\dot x=A(t)x$ and $h(t)$ is continuous.

Then all solutions of 
\begin{equation}\label{eq_perturbed} 
\dot x=A(t)x+h(t)
\end{equation} 
converge to $0$ as $t\to\infty$ if $h\in{\cal{AD}}.$ 

Moreover, if $A(t)$ is entry-wise bounded, then all solutions of (\ref{eq_perturbed}) converge to $0$ as $t\to\infty$ if and only if $h\in{\cal{D}}.$
\end{thm}

\begin{pf}
The statements of theorem are the corollaries of Theorem~A (the part (iv) with $g_2\equiv0$) for general case of $A(t),$ and Theorem~B in \cite{Strauss_Yorke} taking into account that the function $f(t,x)=A(t)x$ with bounded on $[0,\infty)$ matrix function $A(t)$ is a globally Lipschitz function on $\mathbb{R}^n$ in the sense of Section~4 of \cite{Strauss_Yorke} with the Lipschitz constant $L=\max\{\norm{A(t)}, \ t\geq0 \}$.
\end{pf}

Because only the sufficient condition for LTV systems with unbounded $A(t)$ was established, for such systems the origin $x=0$ could be globally attractive also for the perturbations $w\notin\cal{AD},$ even monotonically unbounded ones as it also turns out in Example~\ref{example2}. First, however, we formulate the following theorem on robustness of asymptotically stable LTV systems, which seems to be a completely new result.
\begin{thm}\label{thm:linear}
Let $A(\cdot):\,[0,\infty)\to\mathbb{R}^{n\times n}$ is a continuous matrix function. 

If for some vector norm in $\mathbb{R}^n$ is
\begin{itemize}
\item[1)] $\lim\limits_{t\to\infty}\int\limits_0^t\mu[A(s)]ds=-\infty,$
\end{itemize}
then $x=0\in\mathbb{R}^n$ is globally asymptotically stable equilibrium point of the unperturbed system $\dot x=A(t)x.$
Moreover, if
\begin{itemize}
\item[2)]  $\mu[A(t)]<0$ for all $t$ sufficiently large, and
\item[3)] for all $x\in\mathbb{R}^n$ and all $t\geq0$ is $\norm{w(x,t)}\leq\norm{\tilde{w}(t)}$ with
\[
\lim\limits_{t\to\infty}\big({\norm{\tilde{w}(t)}}/{\mu[A(t)]}\big)=0
\]

\big(or alternatively expressed in the "little-o" notation, 
$
\norm{\tilde{w}(t)}=o\left(\mu[A(t)]\right)\big),
$
\end{itemize}
then all solutions of perturbed system $\dot x=A(t)x+w(x,t)$ converge to $0$ as $t\to\infty.$
\end{thm}
\begin{pf}
Using Lagrange's variation of constants formula for a state-transition matrix  $\Phi_{\mathrm{trans}}(t,\tau)=\Phi(t)\Phi^{-1}(\tau)$ we have that
\begin{equation*}
x(t)=\Phi_{\mathrm{trans}}(t,0)x(0)+\int\limits_{0}^t \Phi_{\mathrm{trans}}(t,\tau)w(x(\tau),\tau)d\tau,
\end{equation*}
and so
\[
\norm{x(t)}\leq\norm{\Phi(t)\Phi^{-1}(0)}\norm{x(0)}+\int\limits_{0}^t \norm{\Phi(t)\Phi^{-1}(\tau)}\norm{w(x(\tau),\tau)}d\tau,
\]
that is, by Lemma~\ref{lognorm_properties} (Item~p3),
\begin{equation*}
\norm{x(t)}\leq \norm{x(0)}e^{\int\limits_0^t\mu[A(s)]ds}+\int\limits_{0}^t e^{\int\limits_{\tau}^t\mu[A(s)]ds}\norm{\tilde w(\tau)}d\tau.
\end{equation*}
We will analyze asymptotics of homogeneous response to initial state and the system's response to external disturbance separately. Obviously, by Assumption~1, 
\[
\lim\limits_{t\to\infty}\norm{x(0)}e^{\int\limits_0^t\mu[A(s)]ds}=0
\]
for an arbitrary initial state $x(0),$  thereby proving the first part of Theorem~\ref{thm:linear} regarding global asymptotic stability of the equilibrium point $x=0$ of unperturbed system ($\tilde w(t)\equiv0$).  Notice that this asymptotic stability is uniform, which is equivalent to uniform exponential stability,  if there exists a real constant $\gamma$ such that $\mu[A(t)]\leq\gamma<0$ for all $t\geq0.$  For the definitions of various types of stability and relations between them, the reader is referred to \cite{Khalil} or, specially for the linear systems, to \cite{Zhou}.

To prove an asymptotic behavior as $t\to\infty$ of the solutions to perturbed system stated in theorem, it remains to analyze the second term on the right-hand side of the above inequality. We have, for $t\geq\tau\geq0,$ that
\[
\int\limits_{0}^t e^{\int\limits_{\tau}^t\mu[A(s)]ds}\norm{\tilde w(\tau)}d\tau= e^{\int\limits_{0}^t\mu[A(s)]ds}\int\limits_{0}^t e^{-\int\limits_{0}^{\tau}\mu[A(s)]ds}\norm{\tilde w(\tau)}d\tau
\]
\[
=\frac{\int\limits_{0}^t e^{-\int\limits_{0}^{\tau}\mu[A(s)]ds}\norm{\tilde w(\tau)}d\tau}{e^{-\int\limits_{0}^t\mu[A(s)]ds}}.
\]
The L'Hospital rule yields 
\[
\lim\limits_{t\to\infty}\frac{\int\limits_{0}^t e^{-\int\limits_{0}^{\tau}\mu[A(s)]ds}\norm{\tilde w(\tau)}d\tau}{e^{-\int\limits_{0}^t\mu[A(s)]ds}}=\lim\limits_{t\to\infty}\frac{\frac{d}{dt}\int\limits_{0}^t e^{-\int\limits_{0}^{\tau}\mu[A(s)]ds}\norm{\tilde w(\tau)}d\tau}{\frac{d}{dt}e^{-\int\limits_{0}^t\mu[A(s)]ds}}
\]
\[
=\lim\limits_{t\to\infty}\frac{e^{-\int\limits_{0}^t\mu[A(s)]ds}\norm{\tilde w(t)}}{e^{-\int\limits_{0}^t\mu[A(s)]ds}(-\mu[A(t)])}=-\lim\limits_{t\to\infty}\frac{\norm{\tilde w(t)}}{\mu[A(t)]},
\]
which together with Assumption~3 gives the claim of Theorem~\ref{thm:linear}.  Observe that if $w(0,t)\equiv0$ for all $t\geq0,$ that is, $x=0$ is an equilibrium point of perturbed system, then under assumptions of Theorem~\ref{thm:linear} $x=0$ is globally asymptotically stable equilibrium point of perturbed system $\dot x=A(t)x+w(x,t).$
\end{pf}
\begin{rmk}
For the LTI systems $\dot x=Ax$ with a Hurwitz matrix $A$, from the end of Example~\ref{example_different_norms}, $\mu_H[A]=-1/\lambda_{\max}(H)$ is negative constant and Assumptions~1 and~2 of Theorem~\ref{thm:linear} are trivially fulfilled. Assumption~3 is satisfied if $\norm{w(x,t)}\leq\norm{\tilde w(t)}\to0$ as $t\to\infty.$ 
\end{rmk}
\begin{ex}\label{example2}
To illustrate Theorem~\ref{thm:linear} by an example, let us consider the system 
\begin{equation}\label{eq:example2}
\dot x=\begin{bmatrix}
-a_1(t) & \beta(t) \\
 -\beta(t) & -a_2(t) 
\end{bmatrix}x
+
\begin{bmatrix}  
t^{7/8} \\
100\cos(t)
\end{bmatrix},\ t\geq0,
\end{equation}
where $a_1(t)=(t+1),$ $a_2(t)=(3+t+\sin t)$ and $\beta(t)$ is an arbitrary continuous function on $[0,\infty).$ 
Using the Euclidean norm, 
\[
\mu_2[A(t)]=\frac12\lambda_{\max}\left(A(t)+A^T(t)\right)=\max\left\{-a_1(t),-a_2(t)\right\}=-(t+1) 
\]
by Table~\ref{table:norms},  

$\int\limits_0^\infty\mu_2[A(s)]ds=-\int\limits_0^\infty (s+1)ds=-\infty$ (Assumption~1 of Theorem~\ref{thm:linear}), 

$\mu_2[A(t)]=-(1+t)<0$ for all $t\geq0$ (Assumption~2), and

$\lim\limits_{t\to\infty}\big({\norm{w(t)}_2}/{\mu_2[A(t)]}\big)=\lim\limits_{t\to\infty}\frac{\left([t^{7/8}]^{2}+[100\cos(t)]^{2}\right)^{1/2}}{-(t+1)}=0$ (Assumption~3).

\centerline{} 

Thus all solutions of the system (\ref{eq:example2}) converge to $0$ for $t\to\infty$ on the basis of Theorem~\ref{thm:linear}. The result of one simulation in the MATLAB environment is shown in Figure~\ref{solution_example2}.

\begin{figure}[ht] 
\captionsetup{singlelinecheck=off}
   \centerline{
    \hbox{
     \psfig{file=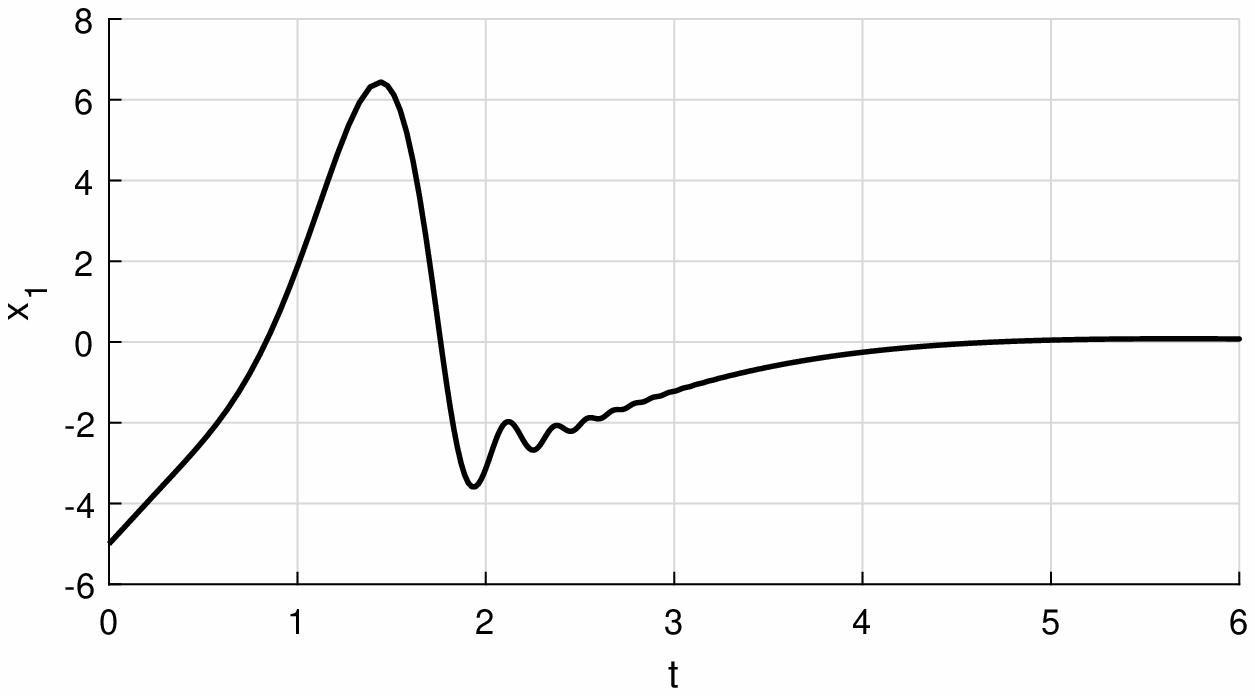,width=5.0cm, clip=}
     \hspace{1.cm}
     \psfig{file=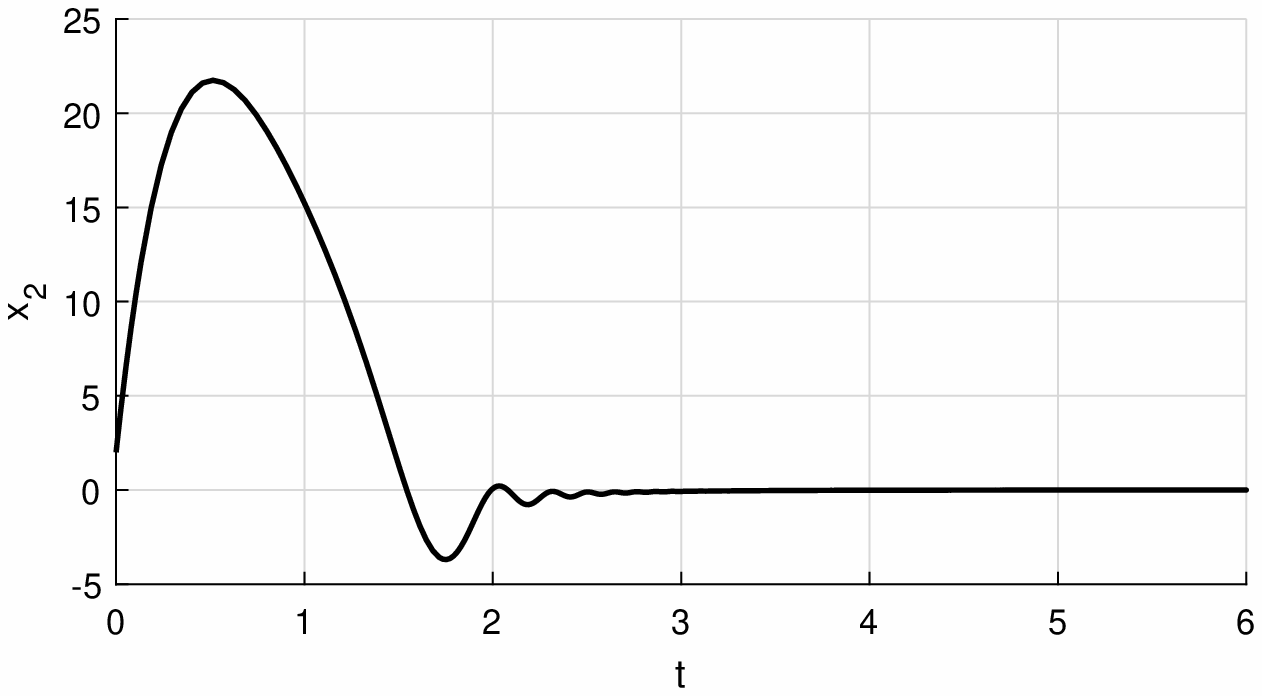,width=5.0cm,clip=}
    }
   }
   
\caption{Solution $x(t)=\left[x_1(t), x_2(t)\right]^T$ of the perturbed system (\ref{eq:example2}) for $\beta(t)=t^4$  and initial state $x(0)=[-5, \ 2]^T.$}
\label{solution_example2}
\end{figure} 

Notice, however, that for $\beta(t)=t^4$ and $a_i(t),$ $i=1,2$ as above, the logarithmic norms $\mu_1[A(t)]=\mu_\infty[A(t)]=t^4-a_1(t),$  and so, the vector norms $\norm{\cdot}_1$ and $\norm{\cdot}_\infty$ are not suitable for stability analysis in this particular case because Assumption~1 of Theorem~\ref{thm:linear} is not satisfied.

We now show that the perturbation $w(t)$ does not belong to the class ${\cal{D}},$ and hence also to the class ${\cal{AD}}.$ For the first component of $w$ and each $\eta>0$ we have 
\[
\int\limits_t^{t+\eta} s^{7/8}ds=\frac{8}{15}\left[(t+\eta)^{15/8}-t^{15/8}\right]=\frac{8}{15}\bigg[\frac{t^{-15/8}-(t+\eta)^{-15/8}}{t^{-15/8}(t+\eta)^{-15/8}}\bigg]
\]
\[ 
=\frac{8}{15}\bigg[\frac{1-\left(1+\eta/t\right)^{-15/8}}{(t+\eta)^{-15/8}}\bigg].
\]
Then the L'Hospital rule yields
\[
\lim\limits_{t\to\infty}\frac{8}{15}\bigg[\frac{1-\left(1+\eta/t\right)^{-15/8}}{(t+\eta)^{-15/8}}\bigg]=\lim\limits_{t\to\infty}\frac{8}{15}\bigg[{\eta}{t^{7/8}}\bigg]\to\infty \ \mathrm{as}\ t\to\infty.
\]
\end{ex} 
\begin{rmk}
In the light of Theorem~\ref{thm:linear} and Example~\ref{example2}, the answer to an open question raised by the authors in \cite{Strauss_Yorke}, below of Theorem~6.1, whether or not $h\in\mathcal{AD}$ is also a necessary condition for perturbations independent of $x$ to converge all solutions of perturbed system (\ref{eq_perturbed}) to $0$ as $t\to\infty,$ is no. Combining Theorems~\ref{LTV} and~\ref{thm:linear}, this question could be modified to the class of functions $h(t)=\tilde w(t)+\tilde h(t),$ where $\tilde w$ satisfies Assumption~3 of Theorem~\ref{thm:linear} and $\tilde h\in\mathcal{AD},$ where the class of allowable perturbations $\tilde w(t)$ could be optimized (maximized) with regard to vector norm in $\mathbb{R}^n.$
\end{rmk}

\begin{ex}
In this example we show, that for perturbation
\[
h(t)=\lambda\left[\sin (e^t),\cos (e^t)\right]^T \in{\cal{D}}\setminus{\cal{AD}},\ \lambda>0,
\]
there exists matrix $A(t)$ (necessarily unbounded, of course) such that $x=0$ is globally (even uniformly) asymptotically stable equilibrium point for $\dot x=A(t)x$ but solution $x^*(t)$ of $\dot x=A(t)x+h(t)$ with $x^*(0)=0$ does not converge to $0$ as $t\to\infty.$ For
\begin{equation*}
A(t)=
\begin{bmatrix} 
-\lambda & e^t\\ 
-e^t & -\lambda 
\end{bmatrix}, \ \lambda>0,
\end{equation*}
obviously, $\mu_2[A(t)]=-\lambda<0$ for all $t\geq0$ proving global uniform asymptotic ($\Leftrightarrow$ uniform exponential) stability of zero solution for unperturbed system, see \cite{Coppel, Khalil, Rugh, Zhou} and Item~p3 of Lemma~\ref{lognorm_properties}. Note that perturbation $h(t)$ does not satisfy Assumption~3 of Theorem~\ref{thm:linear}.

The fundamental matrix for $\dot x=A(t)x$ satisfies
\begin{equation*}
\Phi(t)=e^{-\lambda t}
\begin{bmatrix} 
\sin (e^t) & -\cos (e^t)\\ 
\cos (e^t) &  \sin (e^t)
\end{bmatrix}, \ 
\Phi^{-1}(t)=e^{\lambda t}
\begin{bmatrix} 
\sin (e^t) & \cos (e^t)\\ 
-\cos (e^t) &  \sin (e^t)
\end{bmatrix}
\end{equation*}
and so
\[
x^*(t)=\Phi(t)\int\limits_{0}^t \Phi^{-1}(\tau)h(\tau)d\tau=(1-e^{-\lambda t})
\begin{bmatrix} 
\sin (e^t) \\ 
\cos (e^t) 
\end{bmatrix}
\not\to 0\ \mathrm{as}\ t\to\infty
\]
because 
\[
\norm{x^*(t)}_2=(1-e^{-\lambda t})\to 1. 
\]
This example, together with the previous one, showed that there is no relationship between the admissible perturbations $h(t)$ for LTV systems with bounded and unbounded $A(t),$ respectively. This result also indicates that, in Theorem~\ref{thm:linear}, Assumption~3 regarding the asymptotics of perturbing term cannot be weakened too much and in the last example cannot be achieved better result of the form $\norm{w(t)}=o(t^{\alpha})$ with $\alpha>0$ for the admissible range of perturbations preserving convergence to zero of all solutions of perturbed system as the one we have obtained for the Euclidean norm, $\norm{w(t)}_2=o(1).$ 
\end{ex}
\section*{Conclusions}
In this paper, we derived the novel and relatively easy-to-use criterion for robust global asymptotic stability of the LTV system $\dot x=A(t)x$ being affected by the perturbations (external disturbances).  
Roughly speaking, all solutions of its perturbation $\dot x=A(t)x+w(x,t)$  converge to $0$ if $\norm{w(x,t)}\leq \norm{\tilde{w}(t)}$ and $\big(\norm{\tilde{w}(t)}/\mu[A(t)]\big)\to 0$ as $t\to\infty,$ where by $\mu[A(t)]$ is denoted the logarithmic norm of the system matrix $A(t).$

The fundamental advantage of the approach based on the use of logarithmic norm is the fact that to estimate the norm of state-transition matrix $\Phi(t)\Phi^{-1}(\tau)$ for unperturbed system $\dot x=A(t)x$  we do not need to know the fundamental matrix solution $\Phi(t)$ and all necessary estimates are based purely on the matrix $A(t)$ entries.


\end{document}